\input amstex
\documentstyle{amsppt}
\nologo
\magnification=1200

\topmatter  %XXXXXXXXXXXX

\title Interchanging homotopy limits and colimits in CAT\endtitle
\author by Guillermo Corti\~nas*\endauthor
\thanks (*) Partially supported by ANPC grant BID802/OC-AR-PICT 2260 
and by grants UBACyT EX014 and TW79.
\endthanks
\address Departamento de
Matem\' atica, Ciudad Universitaria Pabell\'on 1, 
(1428) Buenos Aires, Argentina.\endaddress
\email gcorti\@dm.uba.ar\endemail
\leftheadtext{Guillermo Corti\~ nas}
\rightheadtext{Homotopy Limits }

\define\CAT{\operatorname{CAT}}
\define\holi{\operatornamewithlimits{holim}}
\define\hoco{\operatornamewithlimits{hocolim}}
\define\HOM{\operatorname{HOM}}
\abstract Let $I$,$J$ be small categories and $C:I\times J@>>>\CAT$
a functor to the category of small categories. We show that
if $I$ has a final object then 
the canonical map $\hoco_J\holi_IC@>>>\holi_I\hoco_JC$ is a
strong homotopy equivalence.\endabstract
\endtopmatter  %XXXXXXXXXXX

%%%%%%%%%%%%%%%%%%%%%%%%%%%DEFINITIONS%%%%%%%%%%%%%%%%%%%%%%%%%%%%%%%%%

%%%%%%%%%%%%%%%%%%%%%%%%%%%%%%%%%%%%%%%%%%%%%%%%%%%%%%%%%%%%%%%%%%%%%%%
\document
Let $C:I@>>>\bold{CAT}$ be a functor going
from a small category $I$ to the large category of all small categories. 
By the {\it homotopy colimit} of $C$ we mean the Grothendieck
construction ([T]):
$$
\hoco_I C:=\int_IC
$$
and by the {\it homotopy limit}
we mean the pullback:
$$
\CD 
\holi_IC@>>>0\\
@VVV @VVidV\\
\HOM(I,\hoco_IC)@>>\pi_*>\HOM (I,I)\\
\endCD
$$
Here $\HOM$ is the
functor category, $\pi_*$ is induced by the natural projection
$\pi:\hoco_IC@>>>I$, $0$ is the category with only 
one map and $id$
maps the only object of $0$ to the identity functor. A reason for
using  the above definitions is that taking nerves one recovers the
usual homotopy (co)limits for simplicial sets, up to homotopy in the
case of $\hoco$ ([T]) and up to isomorphism in the case of $\holi$ ([L]).
Before we state the main result of this paper, we need a definition.
\bigskip
\definition{Definition} A \it pseudo final\rm\ \ object in a category $I$ 
is a an object $e\in I$ together with a natural map $\epsilon:1@>>>e$ going
from the identity to the constant functor.
\enddefinition
\smallskip
\proclaim{Theorem 1} Let $I$ and $J$ be small categories and 
let $C:I\times J@>>>\CAT$ be a functor. Then there is a 
{\rm faithful} functor:
$$
\iota:\hoco_J\holi_IC@>>>\holi_I\hoco_JC
$$
which is natural in all variables involved.
If in addition the category $I$ has a pseudo-final object, 
then there also exists a functor:
$$
p:\holi_I\hoco_JC@>>>\hoco_J\holi_IC
$$
such that $p\iota=1$ and a natural map $1@>>>\iota p$. The functor
$p$ is natural in $J$ and $C$ and for functors $I@>>>I'$ which preserve
$e$ and $\epsilon$. In particular, under this assumption, the functor $\iota$ 
is a homotopy equivalence.
\endproclaim
\bigskip
Before we prove the theorem we explain our choice of notation and recall
a description of $\holi_IC$ in terms of objects and arrows.
We use a somewhat nonstandard notation for functors $I@>>>\CAT$. On objects
we use subscripts instead of parenthesis; thus $C_i$ (and not $C(i)$)
means the value of $C$ at the object $i\in I$. Also, we omit the letter
$C$ on arrows, so we use the same letter ($\alpha$, $\beta$, $\gamma$,$\dots$)
for a map in $I$ and for its image through the functor $C$. As we use
$i$, $j$, $k$, $\dots$ for objects in $I$ and $x_i$, $y_i$, $z_i$, $\dots$
for objects of $C_i$, this should arise no confusion. For example 
$\alpha:i@>>>j$ is a map in $I$ while $\alpha x_i$ is the object of $C_j$
obtained by applying the image of $\alpha$ through $C$ to the object 
$x_i\in C_i$; i.e. $\alpha x_i$ really means $C(\alpha)(x_i)$. If $\rho:x_i@>>>y_i$
is a map in $C_i$ and $\alpha:i@>>>j$ is in $I$, we write $\alpha(\rho)$ for
$C(\alpha)(\rho)$. 
\bigskip
Here is a description of both $\hoco_IC$ and $\holi_IC$ in terms of
objects
and arrows. An object of $\hoco_IC$ is a pair $x_i:=(i,x)$ where $i\in I$ and
$x\in C_i$. A map $x_i@>>>y_j$ is a pair 
$(\alpha,\rho)$
with $\alpha:i@>>>j\in I$ and  $\rho:\alpha x@>>>y\in C_j$. Composition
is defined as in a semidirect product: 
$(\alpha,\rho)(\beta,\mu)=(\alpha\beta, \alpha(\rho)\mu)$.
On the other hand, $\holi_IC$ is the category of
all pairs of families 
$(x,\rho):=(\{ x_i\}_{i\in I}, \{\rho_\alpha\}_{\alpha\in I})$, indexed 
respectively
 by the objects
and the maps of $I$, where $x_i\in C_i$, and if
$\alpha:i@>>>j$ is a
map in $I$, then $\rho_\alpha:\alpha x_i@>>>x_j$
is a map in
$C_j$. The family $\{\rho_\alpha\}_{\alpha\in I}$ is subject to the conditions:
$$
\rho_1=1
\text{\qquad\qquad}\rho_{\alpha\beta}=\rho_\alpha\alpha(\rho_\beta)\tag1
 $$
In the first
equality, the $1$ on the left is an identity map $1:i@>>>i$
 and the $1$ on the right is the identity of $x_i$ in $C_i$; in the second
equality, $i_0@<\alpha<<i_1@<\beta<<i_2$ are composable maps in $I$
and it should be interpreted as an equality in the set of maps 
$\alpha\beta x_{i_2}@>>>x_{i_0}$ in
$C_{i_0}$. A map $f:(x,\rho^x)@>>>(y,\rho^y)$ in $\holi_IC$ is a
family of maps $f_i:x_i@>>>y_i\in C_i$ indexed by the objects of
$I$ such that the following diagram commutes for every map
$\alpha:i@>>>j\in I$:
$$
\CD
\alpha x_i@>f_i>>\alpha y_i\\
@V\rho^xVV @VV\rho^yV\\
x_j@>>f_j>y_j\\
\endCD\tag2
$$
\bigskip
\demo{Proof of Theorem 1} To begin with, we write down what each of the categories
involved is. An object of $\hoco_I\holi_JC$ is a pair $(j,x)$ where
$j\in J$ is an object and $x=\{(x_{ij},{\rho_\alpha}^x):i,\alpha\in I\}$ is 
a family of objects $x_{ij}\in C_{ij}$, one for each $i\in I$, together
with a family of maps $\rho_\alpha:\alpha x_{ij}@>>>x_{i',j}\in C_{i'j}$,
one for each map $\alpha:i@>>>i'\in I$. (Hereafter we shall write $x=\{(x_{ij},{\rho_\alpha}^x\})$, omitting
the specification of where the indexes lie.) The $\rho_\alpha$ satisfy
condition \thetag{1} for maps in $I$. A map between an object
$(j,x)$ and an object $(k,y)$ is a pair $(\beta, f)$ where 
$\beta:j@>>>k\in J$ and $f=\{f_i\}$ is a family of maps
$f_i\beta x_{ij}@>>>y_{ik}$ such that diagram \thetag{2} commutes. Maps
are composed by $(\gamma,g)(\beta,f)=(\gamma\beta,\{g_i\gamma f_i:i\in I\})$.
An object of the category $\holi_I\hoco_JC$ is a family of pairs 
$(j(i),x_{ij(i)})$ ($i\in I$) where $j(i)\in J$ and
$x_{ij(i)}\in C_{ij(i)}$, together with a family of pairs
$(\beta_\alpha,\rho_\alpha)=(\beta_\alpha^x,\rho_\alpha^x)$, 
($\alpha:i@>>>i'\in I$) where $\beta_\alpha:j(i)@>>>j(i')\in J$ and
$\rho_\alpha:\alpha\beta_\alpha x_{ij(i)}@>>>x_{{i}',j({i}')}\in C_{{i}',j({i}')}$
are maps, and $\rho_1$ and $\beta_1$ are identity maps. The $\beta$
satisfy $\beta_{{\alpha}'\alpha}=\beta_{{\alpha}'}\beta_\alpha$ and
the $\rho$ satisfy \thetag{1}. A map between the object 
$x=((j(i),x_{ij(i)}),(\beta_\alpha,\rho_\alpha^x))$ and the object
$y=((k(i),y_{i,k(i)}),(\beta,\rho_\alpha^y))$ is a family of pairs
$(\beta_i,f_i)$ ($i\in I$), where $\beta_i:j(i)@>>>k(i)\in J$ and
$f_i:\beta_ix_{ij(i)}@>>>y_{ik(i)}\in C_{ik(i)}$ are maps. The
$(\beta_i,f_i)$ satisfy ${\beta_\alpha}^y\beta_i=\beta_{{i}'}{\beta_\alpha}^x$
and $\rho_\alpha^y{\beta_\alpha}^y\alpha(f_i)=f_{{i}'}\beta_{{i}'}(\rho_\alpha^x)$
for each map $\alpha:i@>>>i'\in I$. Composition is defined as
$(\gamma,g)(\beta,f)=\{ (\gamma_i\beta_i,g_i\gamma_i(f_i))\}$. We define
$\iota$ as sending the object $(j,x)$ to the
object $(\{(j,x_{ij})\},\{(1_j,\rho_\alpha^x)\})$ and the map $(\beta,f)$ 
to the map $(\beta,\{f_i\})$. It is clear from the definition that $\iota$
is a functor that maps $\hoco_J\holi_IC$ faithfully into $\holi_I\hoco_JC$.
Moreover the functor $\iota$ identifies $\hoco_J\holi_IC$
with the subcategory of $\holi_I\hoco_JC$ of those objects $x$ as above
which have constant $j=j(i)$ and constant $\beta_\alpha=1_j$, and
with arrows the families $\{(\beta_i,f_i)\}$ with $\beta_i=\beta$, a 
constant
map. Now assume $I$ has a pseudo-final object $\epsilon:1@>>>e$. Define
$p$ as
$p(\{(j(i),x_{ij(i)})\},\{(\beta_\alpha,\rho_\alpha)\})=(\{(j(e),\beta_{{\epsilon}_i}x_{ij(i)})\},\{(1_{j(e)},\beta_{{\epsilon}_i}(\rho_\alpha))\})$
on objects and by $p\{(\beta_i,f_i)\}=(\beta_e,\{\beta_{{\epsilon}_i}(f_i)\})$
on arrows. It is tedious but straightforward to verify that
$p$ is a functor from $\holi_I\hoco_J C$ to 
$\hoco_J\holi_IC$. Once this is verified, it is clear that 
$p\iota=1$. The natural map $\theta:1@>>>\iota p$ is defined as follows.
Given an object $x\in\holi_I\hoco_JC$, define 
$\theta(x)_i=\{(\beta_{{\epsilon}_i}^x,1_{{\beta}_{{\epsilon}_i{x}_{ij(i)}}})\}:\{(j(i),x_{ij(i)})\}@>>>\{(j(e),{\beta_{{\epsilon}_i}}^xx_{ij(i)})\}$.
Another tedious but straightforward verification shows that $\theta$ is
a natural map.\qed
\enddemo
\smallskip
\example{Warning} Write $NC$ for the nerve of $C$. Note that the theorem
does not imply that
$\hoco_J\holi_INC\approx\holi_I\hoco_JNC$. 
This is because, unlike $\holi$,
 $\hoco$ commutes with nerves only up to weak equivalence, not isomorphism,
and $\holi$, unlike $\hoco$, does not preserve all weak equivalences,
only those between fibrant simplicial sets. Theorem 1.2 may
still be true for simplicial sets but it certainly cannot be derived
in this way. Also note that, even when $C$ is fibrant, $\hoco_JC$
need not be so, and thus the spectral sequence for $\holi$ may converge
to the wrong homotopy type. The following example illustrates these
pathologies.
\endexample
\smallskip
\example{Tricky Example} Fix a prime number $p$. Let 
$C:\Bbb N^{op}\times\Bbb N@>>>\CAT$ be given by the following diagram:
$$
\CD
@VVV      @VVV    @VVV\\
(\Bbb Z/p^3)_\delta @>p>>(\Bbb Z/p^3)_\delta @>p>>(\Bbb Z/p^3)_\delta @>p>>\dots\\
@VVV      @VVV    @VVV\\
(\Bbb Z/p^2)_\delta @>p>>(\Bbb Z/p^2)_\delta @>p>>(\Bbb Z/p^2)_\delta @>p>>\dots\\
@VVV      @VVV    @VVV\\
(\Bbb Z/p)_\delta @>p>>(\Bbb Z/p)_\delta @>p>>(\Bbb Z/p)_\delta @>p>>\dots\\
\endCD
$$
Here $(\Bbb Z/p^m)_\delta$ is the discrete category, the vertical maps
are the natural projections, and $p$ means `multiply by $p$'. One checks
that 
$\hoco_{\Bbb N}\holi_{\Bbb N^{op} C\delta}\approx$
\linebreak
$\holi_{\Bbb N^{op}}\hoco_{\Bbb N} (C)_\delta\approx (\Bbb Q_p)_\delta$, the
discrete category of the $p$-adic field. However the $\hoco$ of
each of the rows has the weak homotopy type of a point, by [BK] and [T].
Note this is not in contradiction with the fact that
$\holi$ preserves weak equivalences of fibrant simplicial sets ([BK]) nor the
fact that it preserves adjoint functors ([L]). 
Indeed  
the category $L_m:=\hoco_{\Bbb N}C_{m,-}$ has $\Bbb N\times\Bbb Z/p^m$
as set of objects, and $\hom((n_0,a_m),(n_1,b_m))=\{*\}$ if both 
$n_0\le n_1$ and $p^{{n}_1-{n}_0}a_m=b_m$ and the empty set otherwise.
Thus $NL_m$ is not fibrant, because not every map in $L_m$ is an
isomorphism, and it does not have initial or  final object, i.e.,
the map $L_m@>>>0$ does not have an adjoint.
\endexample

\enddocument